\documentclass{amsart}

\newtheorem{theorem}{Theorem}[section]
\newtheorem{lemma}[theorem]{Lemma}
\newtheorem{proposition}[theorem]{Proposition}
\newtheorem{problem}[theorem]{Problem}

\theoremstyle{definition}
\newtheorem{definition}[theorem]{Definition}
\newtheorem{example}[theorem]{Example}

\theoremstyle{remark}
\newtheorem{remark}[theorem]{Remark}

\numberwithin{equation}{section}
\renewcommand{\dim}{\mathrm{dim}}
\newcommand{\codim}{\mathrm{codim}}

\renewcommand{\span}{\mathrm{span}}
\newcommand{\diam}{\mathrm{diam}} 
\newcommand{\dist}{\mathrm{dist}}
\newcommand{\conv}{\mathrm{conv}}

\newcommand{\R}{\mathbb{R}}
\newcommand{\N}{\mathbb{N}}

\newcommand{\X}{\mathbf{X}}
\renewcommand{\H}{\mathbf{H}}
\newcommand{\Y}{\mathbf{Y}}
\newcommand{\Z}{\mathbf{Z}}
\newcommand{\B}{\mathbf{B}}
\newcommand{\I}{\mathbf{I}}

\newcommand{\aut}{\mathrm{Aut}}
\renewcommand{\S}{\mathbf{S}}

\renewcommand{\ker}{\mathrm{Ker}}
\newcommand{\quotient}{/}
\newcommand{\1}{\boldsymbol{1}}

\begin{document}

\title{On Asymptotic Transitivity in Banach Spaces}

\author{Jarno Talponen}
\address{University of Helsinki, Department of Mathematics and Statistics, Box 68, (Gustaf H\"{a}llstr\"{o}minkatu 2b) FI-00014 University 
of Helsinki, Finland}
\email{talponen@cc.helsinki.fi}

\date{\today}

\begin{abstract}
We introduce a flexible almost isometric version of the almost transitivity property of Banach spaces.
With the help of this new notion we generalize to several directions a strong 
recent rotational characterization of Hilbert spaces due to Randrianantoanina. 
This characterization is a partial answer to the classical Banach-Mazur rotation problem.
\end{abstract}

\maketitle

\section{Introduction}
A Banach space $\X$ is called \emph{transitive} if for each $x\in \S_{\X}$ the orbit
$\mathcal{G}_{\X}(x)=\{T(x)|\ T\colon \X\rightarrow \X\ \mathrm{is\ an\ isometric\ automorphism}\}=\S_{\X}$. 
If $\overline{\mathcal{G}_{\X}(x)}=\S_{\X}$ for all $x\in\S_{\X}$ then $\X$ is called \emph{almost transitive}. 
The following classical \emph{Banach-Mazur rotation problem}, which already appears in Banach's book 
\cite[p.242]{Ba}, remains unsolved despite longstanding active research:
\begin{flushleft}
\textit{Is every separable transitive Banach space $\X$ in fact isometrically a Hilbert space?}
\end{flushleft}

Various partial answers to this problem are known and together with related results they already form a substantial 
theory of rotations in Banach spaces. A recent comprehensive survey of the field is found in \cite{GR}. 
The following recent partial answer to the problem is due to B. Randrianantoanina \cite[Thm.1.1]{Ra2}:
\begin{flushleft}
\textit{If $\X$ is an almost transitive Banach space, which contains a $1$-codimensional $1$-complemented subspace $\Z\subset \X$,
then $\X$ is isometric to a Hilbert space.}
\end{flushleft}

The above result is very general among $1$-codimensional characterizations of Hilbert spaces in terms of rotations, and its 
proof implicitly applies $1$-dimensional linear closest point selections. One should note that the assumptions of the above
result are isometric in their nature. In this paper we show that the rigidity provided by the isometric conditions can be relaxed 
and that the corresponding \emph{almost isometric} conditions are actually sufficient for the characterization 
(see e.g. Theorem \ref{Beata3} below). As an application we also obtain some additional information about the structure of projections 
onto $1$-codimensional and other subspaces.

In section 2 we will introduce and study an asymptotic transitivity property for Banach spaces (see Def. \ref{atdef}) 
which implicitly appears in \cite{Lu}, \cite{Ca2} and \cite{PT} in connection with the \emph{universal disposition property} 
(see \cite{Gu}). In the framework of rotations this \emph{asymptotic transitivity} property seems to be a natural generalization 
of almost transitivity. Thus the asymptotic transitive setting also yields some new information about classical almost transitive spaces.
\subsection{Generalities}
We will use the following notations. Real Banach spaces are denoted by $\X,\Y$ and $\Z$ unless otherwise stated. 
The unit ball and the unit sphere of $\X$ are denoted by $\B_{\X}$ and $\S_{\X}$ respectively. The space of
continuous linear operators $T\colon \X\rightarrow \Y$ is denoted by $L(\X,\Y)$ (abbreviated to $L(\X)$ if $\X=\Y$).
Denote 
\[\mathrm{Aut}(\X)=\{T\in L(\X)|T\ \mathrm{is\ an\ isomorphism}\}.\]
The \emph{Banach-Mazur distance} of mutually isomorphic Banach spaces $\X$ and $\Y$ is given by 
\[d_{BM}(\X,\Y)=\inf\{||T||\cdot ||T^{-1}||:\ T\in L(\X,\Y)\ \mathrm{isomorphism}\}.\]
Spaces $\X$ and $\Y$ are said to be \emph{almost isometric} if $d_{BM}(\X,\Y)=1$.  
Denote the group of rotations of $\X$ by 
\[\mathcal{G}_{\X}=\{T\in\mathrm{Aut}(\X)|\max(||T||,||T^{-1}||)=1\}\] 
and the \emph{orbit} of $x\in\S_{\X}$ by $\mathcal{G}_{\X}(x)=\{T(x)|T\in\mathcal{G}_{\X}\}$. 
If $\overline{\conv}(\mathcal{G}_{\X}(x))=\B_{\X}$ for all $x\in\S_{\X}$ then $\X$ is called \emph{convex-transitive}.
If for any equivalent norm $|||\cdot|||\sim ||\cdot||$ such that $\mathcal{G}_{(\X,||\cdot||)}\subset\mathcal{G}_{(\X,|||\cdot|||)}$
it holds that $\mathcal{G}_{(\X,||\cdot||)}=\mathcal{G}_{(\X,|||\cdot|||)}$ then $||\cdot||$ is called a \emph{maximal} norm.
We denote by $\I=\I_{\X}$ the identical mapping.
In a metric space $(X,d)$ the Hausdorff distance of non-empty sets
$A,B\subset X$ is defined as 
\[d_{H}(A,B)=\sup\{\dist(a,B),\dist(b,A)|\ a\in A, b\in B\}.\]
We say that $\Y,\Z\subset \X$ are $L^{p}$-summands of $\X$ if $\X=\Y\oplus\Z$, where 
$||(y,z)||=(||y||_{_{\Y}}^{p}+||z||_{_{\Z}}^{p})^{\frac{1}{p}}$ (respectively $\max(||y||_{\Y},||z||_{\Z})$ for $p=\infty$)  
for all $(y,z)\in\Y\times\Z$. We denote this by $\X=\Y\oplus_{p}\Z$. 
We write shortly $L^{p}=L^{p}((0,1),m)$, where $m$ is the Lebesgue measure on $(0,1)$. 
For $f\in \X^{\ast}$ and $x\in\X$ we let $f\otimes x\colon \X\rightarrow [x]$
be the map $y\mapsto f(y)x$. 

For a discussion of basic concepts and results concerning the geometry of the norm we refer to the first chapter of \cite{JL}.
Recall that a norm $||\cdot||$ is said to be \emph{Gateaux differentiable} at $x\in \S_{\X}$ if there is
$l\in \X^{\ast}$ such that 
$$\lim_{t\rightarrow 0}\frac{||x+th||-1}{t}=l(h)$$
for every $h\in \X$. If $||\cdot||$ is Gateaux differentiable at $x\in \S_{\X}$ and the derivative above satisfies
$$\lim_{||h||\rightarrow 0}\frac{||x+h||-1}{||h||}-l(h)=0$$
then $||\cdot||$ is said to be Fr\'{e}chet differentiable at $x\in\S_{\X}$. 
Alternatively, the point $x$ is called above \emph{(Gateaux)-smooth}, respectively \emph{Frechet-smooth}.
Let us recall the following two classical results (see e.g. \cite[p.92,300]{HHZ}).
\begin{lemma}(Smulyan)\label{lm FS}
The following conditions $(1)-(2)$ are equivalent:
\begin{enumerate}
\item[(1)]{The norm $||\cdot||$ of a Banach space $\X$ is Frechet differentiable at $x\in \S_{\X}$.}
\item[(2)]{For all $(f_{n}),(g_{n})\subset \S_{\X^{\ast}}$ such that $\lim_{n\rightarrow\infty}f_{n}(x)=\lim_{n\rightarrow\infty}g_{n}(x)=1$ it holds that $\lim_{n\rightarrow\infty}||f_{n}-g_{n}||=0$.}
\end{enumerate}
\end{lemma}
\begin{lemma}(Smulyan)\label{lm GS}
The following conditions $(1)-(3)$ are equivalent:
\begin{enumerate}
\item[(1)]{The norm $||\cdot||$ of a Banach space $\X$ is Gateaux differentiable at $x\in \S_{\X}$.}
\item[(2)]{For all $(f_{n}),(g_{n})\subset \S_{\X^{\ast}}$ such that 
$\lim_{n\rightarrow\infty}f_{n}(x)=\lim_{n\rightarrow\infty}g_{n}(x)=1$ it holds that $f_{n}-g_{n}\stackrel{\omega^{\ast}}{\longrightarrow} 0\ \mathrm{as}\ n\rightarrow\infty$.}
\item[(3)]{There is unique $f\in \S_{\X^{\ast}}$ such that $f(x)=1$.}
\end{enumerate}
\end{lemma}
\begin{theorem}(Bishop-Phelps)\label{BPT}
The norm-attaining functionals of $\S_{\X^{\ast}}$ are dense in $\S_{\X^{\ast}}$.
\end{theorem}

\section{Asymptotic transitivity and structures of orthogonal type}

\begin{definition}\label{atdef}
We say that a Banach space $\X$ is \emph{asymptotically transitive} 
if for all $x\in \S_{\X}$ the generalized orbit $\mathcal{O}_{\X}(x)$ of $x$ defined by
\[\mathcal{O}(x)=\mathcal{O}_{\X}(x)=\bigcap_{\epsilon>0}\{T(x)|T\in \mathrm{Aut}(\X),\ \max(||T||,||T^{-1}||)\leq 1+\epsilon\}\]
satisfies $\mathcal{O}_{\X}(x)=\S_{\X}$.
\end{definition}

We will first list some basic facts about asymptotically transitive spaces:
\begin{remark}\label{remark1}\

\begin{enumerate}
\item[(i)]{$\mathcal{O}(x)$ is norm closed for all $x\in\S_{\X}$ (see also the proof of Prop. \ref{th eta}),}
\item[(ii)]{almost transitive spaces are asymptotically transitive (see also \cite[p. 239]{Lu}).}
\item[(iii)]{If $\X$ and $\Y$ are almost isometric and $\X$ is almost transitive, then $\Y$ is asymptotically transitive.}
\item[(iv)]{The convex-transitive space $L^{\infty}(0,1)$ is not asymptotically transitive.}
\end{enumerate}
\end{remark}
The above claims (i)-(iii) are fairly immediate, where (i) implies (ii). For the last claim, observe that if 
$T\in \mathrm{Aut}(L^{\infty})$
then $\mathrm{ess\ inf}_{t}|T(x)(t)|\geq ||T^{-1}||^{-1} \mathrm{ess\ inf}_{t}|x(t)|$ for 
all $x\in L^{\infty}(0,1)$.
Hence $\chi_{[0,\frac{1}{2}]}\notin \mathcal{O}_{L^{\infty}}(\chi_{[0,1]})$ even though $L^{\infty}$ is
convex-transitive (see \cite[p.17]{GR}).

The generalized orbit $\mathcal{O}(x)$ and the closure $\overline{\mathcal{G}(x)}$ need not coincide in general 
(see Example \ref{example XYZ}), but unfortunately we do not know at the moment any examples of asymptotically transitive spaces, 
which are not almost transitive.

Our work in this section is aimed at showing the following main result, which we will prove after some auxiliary results.
\begin{theorem}\label{Beata3}
Let $\X$ be an asymptotically transitive Banach space such that for each $\epsilon>0$ there is $1$-codimensional
$1+\epsilon$-complemented subspace $Z_{\epsilon}\subset\X$. Then $\X$ is in fact isometrically a Hilbert space.
\end{theorem}

\newcommand{\ei}{\overset{_{\wedge}}{\boldsymbol{\eta}}}
\newcommand{\es}{\overset{_{\vee}}{\boldsymbol{\eta}}}
We will first establish some crucial observations. Towards Theorem \ref{Beata3} define the local projectional indices 
$\ei_{ _{\X}},\es_{ _{\X}}\colon \S_{\X}\rightarrow \R$ for $y\in\S_{\X}$ by
\[\ei(y)=\ei_{ _{\X}}(y)=\inf\{||\I-f\otimes y||:\ f\in\X^{\ast},\ f(y)=1\},\]
\[\es(y)=\es_{ _{\X}}(y)=\inf_{\epsilon>0}\sup\{||\I-f\otimes y||:\ f\in (1+\epsilon)\B_{\X^{\ast}},\ f(y)=1\}.\]
Clearly $1\leq\ei_{ _{\X}}\leq\es_{ _{\X}}\leq 2$ and 
$\ei_{ _{\Y}}\leq\ei_{ _{\X}}$, $\es_{ _{\Y}}\leq\es_{ _{\X}}$ pointwise for $\Y\subset\X$.

\begin{proposition}\label{th eta}
For any Banach space $\X$ the following holds:
\begin{enumerate}
\item[(i)]{The map $\ei$ is uniformly continuous and $\es$ is upper semicontinuous on $\S_{\X}$.}
\item[(ii)]{For each $x\in \S_{\X}$ there is $f\in\X^{\ast}$ such that $f(x)=1$ and $||\I-f\otimes x||=\ei(x)$.}
\end{enumerate}
\end{proposition}
\begin{proof}
First we check claim (ii). Let $x\in\S_{\X}$ and $(f_{n})\subset \X^{\ast}$ be a sequence such that $f_{n}(x)=1$ and 
\mbox{$||\I-f_{n}\otimes x||\leq \ei(x)+n^{-1}$} for each $n\in\N$. Note that 
\[||f_{n}||=||f_{n}\otimes x||\leq 1+||\I-f_{n}\otimes x||\leq 1+\ei(x)+n^{-1},\]
so that $(f_{n})\subset 5\B_{\X^{\ast}}$. By the Banach-Alaoglu theorem $(5\B_{\X^{\ast}},\omega^{\ast})$ is compact, and thus 
there exists a $\omega^{\ast}$-cluster point $f\in 5\B_{\X^{\ast}}$ for the sequence $(f_{n})\subset 5\B_{\X^{\ast}}$. Clearly $f(x)=1$ and
it is not difficult to check that $||\I-f\otimes x||\leq \limsup_{n\rightarrow\infty}||\I-f_{n}\otimes x||=\ei(x)$, so that the proof of claim (ii) is complete.

Next we will make some general observations. Since $1\leq \ei \leq \es\leq 2$, we may assume without loss of
generality as above, that each rank-$1$ map $h\otimes z$, where $h\in \X^{\ast},\ z\in\S_{\X},$ appearing in this proof 
satisfies $||h\otimes z||\leq 5$. Let $0<\epsilon<\frac{1}{5}$ and $x,y\in\S_{\X}$ be such that $||x-y||=\epsilon$.
Let $f\in \X^{\ast}$ be such that $f(x)=1$ and $||f||\leq 5$. Define $S\in L(\X)$ by $S=\I+f\otimes(y-x)$. 
Note that $S(x)=y$. Since $||x-y||=\epsilon<\frac{1}{5}$ we obtain that $|f(x)-f(y)|\leq ||f||\cdot ||x-y||<1$. Thus
$f(y)>0$. Hence $\ker(f)\oplus [x]=\ker(f)\oplus [y]=\X$ as $x,y\notin\ker(f)$. 
Clearly $S_{|\ker(f)}=\I_{|\ker(f)}$. Thus $S\colon \ker(f)\oplus [x]\rightarrow \ker(f) \oplus [y]$ is a bijection.
Note that $||v-S(v)||\leq |f(v)|\cdot ||x-y||$. Thus
\begin{eqnarray} 
||S||\leq ||\I||+||\I-S||\leq 1+||f||\cdot ||x-y||\leq 1+5\epsilon\label{eq: S5eps}\\
||S^{-1}||\leq (1-||f||\cdot ||x-y||)^{-1}\leq (1-5\epsilon)^{-1},\label{eq: SS-1}
\end{eqnarray}
since $||Sv||\geq ||v||-||v-Sv||\geq (1-||f||\cdot ||x-y||)||v||$ for $v\in\X$.
Hence
\begin{equation}\label{eq: SF}
\max(||S||,||S^{-1}||)\leq \delta(\epsilon),
\end{equation}
where $\delta(\epsilon)\stackrel{\cdot}{=}\max(1+5\epsilon,(1-5\epsilon)^{-1})$ for $0<\epsilon<\frac{1}{5}$.
From \eqref{eq: S5eps} and \eqref{eq: SS-1} it follows that
\begin{equation}\label{eq: I-S-1}
||\I-S^{-1}||=||(\I-S)S^{-1}||\leq ||\I-S||\cdot ||S^{-1}|| \leq  \frac{5\epsilon}{1-5\epsilon}.
\end{equation}
Note that since $P=f\otimes x\colon \X\rightarrow [x]$ is a linear projection then so is
$S\circ P\circ S^{-1}\colon \X\rightarrow [y]$. Moreover, 
\begin{equation}\label{eq: ISPS}
||S\circ S^{-1}-S\circ P\circ S^{-1}||\leq||S||\cdot ||\I-P||\cdot ||S^{-1}||\leq \delta(\epsilon)^{2}||\I-P||.
\end{equation}
Note that above $\delta(\epsilon)\rightarrow 1$ as $\epsilon\rightarrow 0^{+}$, so that by (\ref{eq: ISPS}) the map 
$\ei$ is uniformly continuous.

For the case $\es$ we observe that by (\ref{eq: SF}) and the fact that $S^{-1}(y)=x$ it holds for all $\alpha\in [0,1]$ that 
\begin{equation}\label{eq: S*}
(S^{-1})^{\ast}(\{f\in (1+\alpha)\B_{\X^{\ast}}| f(x)=1\})\subset \{g\in (1+\alpha)\delta(\epsilon)\B_{\X^{\ast}}|g(y)=1\}.
\end{equation}
Hence by applying (\ref{eq: SF}), (\ref{eq: S*}), \eqref{eq: S5eps} for $g$ and (\ref{eq: I-S-1}) respectively, we obtain that
\begin{equation*}
\begin{array}{l}
\sup\{||\I-f\otimes x||:\ f\in (1+\alpha)\B_{\X^{\ast}},\ f(x)=1\}\\
\leq\delta(\epsilon)\sup\{||(\I-f\otimes x)S^{-1}||:\ f\in (1+\alpha)\B_{\X^{\ast}},\ f(x)=1\}\\
\leq \delta(\epsilon)\sup\{||S^{-1}-g\otimes x||:g\in (1+\alpha)\delta(\epsilon)\B_{\X^{\ast}}, g(y)=1\}\\
\leq \delta(\epsilon)\sup\{||S^{-1}-g\otimes y||:g\in (1+\alpha)\delta(\epsilon)\B_{\X^{\ast}}, g(y)=1\}+(1+\alpha)\delta(\epsilon)\epsilon\\
\leq \delta(\epsilon)\sup\{||\I-g\otimes y||:g\in (1+\alpha)\delta(\epsilon)\B_{\X^{\ast}}, g(y)=1\}+(1+\alpha+\frac{5}{1-5\epsilon})\delta(\epsilon)\epsilon.\\
\end{array}
\end{equation*}
Let $z\in\S_{\X}$ and $(z_{n})\subset \S_{\X}$ be a sequence such that $||z-z_{n}||=\beta_{n}\rightarrow 0$ as $n\rightarrow\infty$. 
By the preceding estimate we obtain that 
\begin{equation}
\begin{array}{l}
\sup\{||\I-g_{n}\otimes z_{n}||:\ g_{n}\in (1+\beta_{n})\B_{\X^{\ast}},\ g_{n}(z_{n})=1\}\\
\leq \delta(\beta_{n})\sup\{||\I-g\otimes z||:g\in (1+\beta_{n})\delta(\beta_{n})\B_{\X^{\ast}},\ g(z)=1\}\\
+(1+\beta_{n}+\frac{5}{1-5\beta_{n}})\delta(\beta_{n})\beta_{n}\\
\rightarrow \es(z)\ \mathrm{as}\ n\rightarrow\infty,
\end{array}
\end{equation}
so that $\limsup_{n\rightarrow\infty}\es(z_{n})\leq \es(z)$, that is, the upper semicontinuity of $\es$ on $\S_{\X}$.
\end{proof}

\begin{proposition}\label{asyconstant}
Suppose that $\X$ is an asymptotically transitive Banach space. Then $\ei$ and $\es$ are constant functions.
\end{proposition}
\begin{proof}
Suppose that $x,y\in\S_{\X}$. By the asymptotic transitivity of $\X$ there is 
a sequence $(T_{n})\subset\mathrm{Aut}(\X)$ such that $T_{n}(x)=y$ and $\max(||T_{n}||,||T_{n}^{-1}||)\leq 1+n^{-1}$
for each $n\in\N$. 
By part (ii) of Proposition \ref{th eta} there is $f\in \X^{\ast},\ f(x)=1$ 
such that $||\I-f\otimes x||=\ei(x)$.
Define $g_{n}=f\circ T_{n}^{-1}$ for $n\in\N$. Then $g_{n}(y)=1$ for $n\in\N$, so that
\[\ei(y)\leq\limsup_{n\rightarrow\infty}||\I-g_{n}\otimes y||\leq\limsup_{n\rightarrow\infty}(||T_{n}||\cdot ||\I-f\otimes x||\cdot||T_{n}^{-1}||)=\ei(x).\] 
By symmetry we obtain that $\ei(x)=\ei(y)$.

In the case of $\es$ suppose similarly as above that $(f_{n})\subset 2\B_{\X^{\ast}}$ is a sequence so that $f_{n}(x)=1$ for 
each $n\in\N$, $||f_{n}||\rightarrow 1$ and $||\I-f_{n}\otimes x||\rightarrow \es(x)$ as $n\rightarrow\infty$.
Then putting $g_{n}=f_{n}\circ T_{n}^{-1},\ n\in\N,$ as above yields the claim. 
\end{proof}
\begin{proof}[Proof of Theorem \ref{Beata3}]
In the case $\dim(\X)\leq 2$ the claim is clear. For the case $\dim(\X)\geq 3$ 
recall that a combination of \cite[Lemma 13.1]{Am} and \cite[Criterion 13.4']{Am} yields the following characterization 
of Hilbert spaces:\\
\textit{A Banach space $E,\ \dim(E)\geq 3,$ is isometrically a Hilbert space if and only if for all 
$e\in\S_{E}$ there is $f\in E^{\ast}$ such that $f(e)=1$ and $||\I-f\otimes e||=1$.}    
\newline
By combining Propositions \ref{th eta} and \ref{asyconstant} we obtain the above condition in our case. 
\end{proof}
One should note that the above result was anticipated (in the almost transitive setting) in F. Cabello's unpublished Ph.D. thesis 
\cite{Ca3}. See also \cite[Prop. 1.5]{AF} for a related result.

\begin{proposition}\label{IP} 
Suppose that $\X$ is asymptotically transitive and Frechet-smooth. Then $\es=C$, a constant function, 
and $||\I-P||=C$ for any $1$-dimensional linear projection $P\in L(\X)$ such that $||P||=1$. 
\end{proposition}
\begin{proof}
Let $y\in\S_{\X}$. 
Suppose that the sequence $(g_{n})\subset 2\B_{\X^{\ast}}$ is chosen as in the proof of Proposition \ref{asyconstant} for $\es$. 
Then an application of the Frechet-smoothness together with the Smulyan lemma gives that
$g_{n}\stackrel{||\cdot||}{\longrightarrow}g$ as $n\rightarrow\infty$ for the unique $g\in \S_{\X^{\ast}}$ such that $g(y)=1$
and $g\otimes y=P$. We obtain that $||\I-g_{n}\otimes y||\rightarrow ||\I-g\otimes y||$ as $n\rightarrow\infty$. 
Thus $||\I-g\otimes y||=\es(y)$. This yields the claim, since $\es$ is a constant function
by Proposition \ref{asyconstant}. 
\end{proof}
It turns out (see Theorem \ref{ucs} below) that in the asymptotically transitive setting rather mild geometric conditions, for example
the Frechet-differentiability of the norm at \emph{some} point $x\in\S_{\X}$, already imply the Frechet-differentiability 
of the norm at \emph{all} $y\in\S_{\X}$.

\section{Asymptotic transitivity and the geometry of the norm}
We will show that in the case of asymptotically transitive spaces some suitable control of the convexity or of the smoothness 
of the unit ball actually quarantees simultaneously \emph{uniform convexity and uniform smoothness}. 
We will use some concepts and results about the geometry of the norm, for a discussion of which we refer to \cite[Ch.1]{JL}.
Denote \emph{closed slices} by $S(\B_{\X},f,\alpha)=\{x\in\B_{\X}|f(x)\geq \alpha\}$ for $f\in\S_{\X^{\ast}},\ 0\leq\alpha\leq 1$.
It is said that $B_{\X}$ is \emph{dentable} (respectively $B_{\X^{\ast}}$ is \emph{$\omega^{\ast}$-dentable}) if 
\[\inf_{\substack{f\in\S_{\X^{\ast}},\\ 0<\alpha<1}}\diam(S(\B_{\X},f,\alpha))=0\ \ (resp. \inf_{\substack{x\in\S_{\X},\\ 0<\alpha<1}}\diam(S(\B_{\X^{\ast}},x,\alpha))=0).\] 
The point $x\in \S_{\X}$ is strongly exposed by $f\in\S_{\X^{\ast}}$ if $f(x)=1$ and\\ 
$\inf_{\alpha<1}\diam(S(\B_{\X},f,\alpha))=0$. Recall that the local modulus of convexity\\ 
$\Delta_{\X}\colon \S_{\X}\times [0,1)\rightarrow \R$ is defined by
\[\Delta_{\X}(z,\epsilon)=\inf\{1-\lambda\in [0,\infty) |\ \mathrm{there\ is}\ v\in \X\ \mathrm{such\ that}\ ||v||\geq \epsilon,\ ||\lambda z\pm v||\leq 1\}\]
for $z\in \S_{\X}$ and $0\leq\epsilon<1$. Note that $\Delta_{\X}(z,0)=0$ for each $z\in \S_{\X}$.
\begin{theorem}\label{ucs}
Let $\X$ be an asymptotically transitive Banach space. If $\B_{\X}$ is dentable or $\B_{\X^{\ast}}$ is $\omega^{\ast}$-dentable, then
$\X^{\ast}$ is also asymptotically transitive, and both $\X$ and $\X^{\ast}$ are uniformly convex and uniformly smooth.
\end{theorem}

\begin{proof}
We will first prove the following auxiliary technical claim:\\ 
\textit{Claim. Suppose that $\X$ and $\Y$ are Banach spaces and $x\in\S_{\X},\ y\in\S_{\Y}$. 
If $T\colon \X\rightarrow\Y$ is a linear isomorphism such that $T(x)=y$, then
\begin{equation}
\Delta_{\Y}\Big(y,\frac{\epsilon\cdot \Delta_{\X}(x,\epsilon)}{||T^{-1}||(||T||+\Delta_{\X}(x,\epsilon)-1)}\Big)\leq \Delta_{\X}(x,\epsilon). 
\end{equation}} 
We apply the convention $\frac{0}{0}=0$ in the fraction in the left side of the formula above. 
Indeed, to verify the claim, let us use the same notations as in the definition of $\Delta_{\X}$. 
Observe that $\min(||T||,||T^{-1}||)\geq 1$.
The claim holds trivially if $\epsilon\cdot\Delta_{\X}(x,\epsilon)=0$.
Suppose that $0<\epsilon<1$, $\Delta_{\X}(x,\epsilon)>0$ and write $1-\lambda=\Delta_{\X}(x,\epsilon)$. Observe that $0\leq \lambda\leq 1$.
Let $(v_{n})$ be a sequence such that 
$||\lambda x\pm v_{n}||\leq 1,\ n\in\N,$ and $\liminf_{n\rightarrow\infty} ||v_{n}||_{\X}\geq \epsilon$. 

Define an auxiliary mapping $\phi\colon \R \rightarrow (0,\infty)$ by 
$\phi(t)=\frac{1-\lambda}{t-\lambda}$ if $t>1$ and $\phi(t)=1$ otherwise. 
Write $T(v_{n})=w_{n},\ n\in\N$. Fix $n\in\N$.

We aim to prove the following subclaim:
\begin{equation}\label{eq: Deltaphi}
\Delta_{\Y}(y,\alpha_{n}||w_{n}||)\leq \Delta_{\X}(x,\epsilon),
\end{equation}
where $\alpha_{n}=\min_{\theta=\pm 1}\phi(||\lambda y+\theta w_{n}||)$.
Indeed, since $h\mapsto \frac{||\lambda y+ hw||-\lambda}{h}$ is a non-decreasing map for each $w$, we obtain by substituting 
$h=\frac{1-\lambda}{||\lambda y+w_{n}||-\lambda}$ once\\ 
$||\lambda y +w_{n}||> 1$ that
\[\frac{\Big|\Big| \lambda y+\frac{1-\lambda}{||\lambda y +w_{n}||-\lambda}w_{n}\Big|\Big| -\lambda}{\frac{1-\lambda}{||\lambda y +w_{n}||-\lambda}}\leq ||\lambda y +w_{n}||-\lambda.\]
Consequently
\[\frac{\Big|\Big| \lambda y+\frac{1-\lambda}{||\lambda y +w_{n}||-\lambda}w_{n}\Big|\Big| -\lambda}{1-\lambda}\leq 1,\ \
\mathrm{and\ hence}\ \ \Big|\Big| \lambda y+\frac{1-\lambda}{||\lambda y +w_{n}||-\lambda}w_{n}\Big|\Big|\leq 1.\]
Hence we obtain for $w_{n}$ in general that 
\[|| \lambda y+\phi(||\lambda y +w_{n}||) w_{n}||\leq 1\] 
and similarly 
\[|| \lambda y-\phi(||\lambda y -w_{n}||)w_{n}||\leq 1.\]
Observe that $\{h\in \R:\ ||\lambda y+h w_{n}||\leq 1\}\subset \R$ is a compact interval that contains $0$ 
by the basic properties of $||\cdot||_{\Y}$. Moreover, the previous estimates give that
\[[-\phi(||\lambda y -w_{n}||),\phi(||\lambda y +w_{n}||)] \subset \{h\in \R :\ ||\lambda y+h w_{n}||\leq 1\}.\]
In particular we obtain that 
\[[-\alpha_{n},\alpha_{n}]\subset \{h\in\R :\ ||\lambda y+h w_{n}||\leq 1\}\]
and 
\[||\lambda y\pm \alpha_{n} w_{n}||\leq 1.\] 
Hence we have proved the subclaim \eqref{eq: Deltaphi}.

Observe next that
\[||\lambda y\pm w_{n}||-\lambda=||T(\lambda x\pm v_{n})||-\lambda \leq ||T||-\lambda=||T||+\Delta_{\X}(x,\epsilon)-1\] 
since by assumption $||\lambda x\pm v_{n}||\leq 1$ and $\Delta_{\X}(x,\epsilon)=1-\lambda$. Thus  
\begin{equation}\label{eq: TD}
\frac{\Delta_{\X}(x,\epsilon)}{||T||+\Delta_{\X}(x,\epsilon)-1}\leq\min_{\theta=\pm 1}\Big(\frac{1-\lambda}{||\lambda y+\theta w_{n}||-\lambda}\Big)=\alpha_{n}.
\end{equation}
Note that above $||T||-1\geq 0$ and hence $\frac{\Delta_{\X}(x,\epsilon)}{||T||+\Delta_{\X}(x,\epsilon)-1}\leq 1$. 
Since $||w_{n}||\geq \frac{||v_{n}||}{||T^{-1}||}$, where $\liminf_{n\rightarrow\infty}\frac{||v_{n}||}{||T^{-1}||}\geq 
\frac{\epsilon}{||T^{-1}||}$ and $\Delta_{\Y}(y,\cdot)$ is non-decreasing by definition, we obtain 
by combining \eqref{eq: Deltaphi} and \eqref{eq: TD} that
\[\Delta_{\Y}\Big(y,\frac{\epsilon\cdot \Delta_{\X}(x,\epsilon)}{||T^{-1}||(||T||+\Delta_{\X}(x,\epsilon)-1)}\Big)\leq \Delta_{\X}(x,\epsilon),\]
which is the Claim. 
 
Let $x,y\in \S_{\X}$. Then $\Delta_{\X}(x,\cdot),\Delta_{\X}(y,\cdot)\colon (0,1)\rightarrow \R$ 
are non-decreasing maps. Thus they are continuous a.e. with respect to the Lebesgue measure $m$ on $(0,1)$ and 
\[m(\{t\in (0,1)|\ \mathrm{both}\ \Delta_{\X}(x,\cdot)\ \mathrm{and}\ \Delta_{\X}(y,\cdot)\ \mathrm{are\ continuous\ at}\ t\})=1.\]
Let $t_{0}\in (0,1)$ be such a point of joint continuity. Then by the asymptotic transitivity of $\X$ there is a sequence 
$(T_{n})\subset \mathrm{Aut}(\X)$ such that $T_{n}(x)=y,\ n\in\N,$ and $||T_{n}||\cdot ||T_{n}^{-1}||\rightarrow 1$ as 
$n\rightarrow\infty$. It follows from the claim and symmetry that $\Delta_{\X}(x,t_{0})=\Delta_{\X}(y,t_{0})$. 
We conclude that $\Delta_{\X}(x,\cdot)$ and $\Delta_{\X}(y,\cdot)$ coincide $m$-a.e. Hence 
\begin{equation}\label{eq: limxy}
\lim_{t\rightarrow s^{-}}\Delta_{\X}(x,t)=\lim_{t\rightarrow s^{-}}\Delta_{\X}(y,t)
\end{equation}
for all $s\in (0,1)$ and $(x,y)\in \S_{\X}\times \S_{\X}$.

Consider first the case where $\B_{\X}$ is dentable. Fix sequences $(f_{n})\subset \S_{X^{\ast}}$ and $(\alpha_{n})\subset (0,1)$ 
such that $\alpha_{n}\rightarrow 1$ and $\diam(S(\B_{\X},f_{n},\alpha_{n}))\rightarrow 0$ as $n\rightarrow\infty$. The norm-attaining 
functionals of $\S_{\X^{\ast}}$ are dense subset by the Bishop-Phelps Theorem (see Theorem \ref{BPT}), so that 
we can find norm-attaining functionals $(g_{n})\subset \S_{\X^{\ast}}$ such that $||f_{n}-g_{n}||<\frac{1-\alpha_{n}}{2}$ for $n\in\N$.
Then $S(\B_{\X},g_{n},\alpha_{n}+\frac{1-\alpha_{n}}{2})\subset S(\B_{\X},f_{n},\alpha_{n})$ for $n\in\N$.
Let $(x_{n})\subset \S_{\X}$ be such that $g_{n}(x_{n})=1$ for $n\in\N$.

The fact that 
\begin{equation}\label{eq: dxe}
\Delta_{\X}(x_{n},\epsilon)\geq \inf\{1-\lambda| \diam(S(\B_{\X},g_{n},\lambda))\geq \epsilon\}\ \mathrm{for\ all}\ 0<\epsilon<1
\end{equation}
is a consequence of the following observation: if $\lambda x_{n}\pm v\in \B_{\X},\ ||v||\geq\epsilon,$ then at least one of 
$\lambda x_{n}\pm v$ is contained in $S(\B_{\X},g_{n},\lambda)$, so that $\diam(S(\B_{\X},g_{n},\lambda))\geq\epsilon$. 
Hence, by choosing $\epsilon=2\diam\big(S\big(\B_{\X},g_{n},\alpha_{n}+\frac{1-\alpha_{n}}{2}\big)\big)$ in \eqref{eq: dxe} we obtain that 
\begin{equation*}
\begin{array}{ll}
\Delta_{\X} \left( x_{n},2\diam \left( S \left( \B_{\X},g_{n},\alpha_{n}+\frac{1-\alpha_{n}}{2} \right) \right) \right)\\
\geq\inf \left\{ 1-\lambda|\ \diam(S(\B_{\X},g_{n},\lambda))\geq 2\diam \left( S \left( \B_{\X},g_{n},\alpha_{n}+\frac{1-\alpha_{n}}{2} \right) \right) \right\}\\
>1-\left( \alpha_{n}+\frac{1-\alpha_{n}}{2} \right) ,\quad n\in\N,\\
\end{array}
\end{equation*}
\ \medskip \ \\
where we used the fact that $\diam(S(\B_{\X},g_{n},\lambda))\geq 2\diam\left(S\left(\B_{\X},g_{n},\alpha_{n}+\frac{1-\alpha_{n}}{2}\right)\right)$ implies
that $\lambda<\alpha_{n}+\frac{1-\alpha_{n}}{2}$. Let us summarize the facts established so far: 
$\Delta_{\X}(x,\cdot)=\Delta_{\X}(y,\cdot)$ holds $m$-a.e. for each pair $(x,y)\in \S_{\X}\times \S_{\X}$, these maps are non-decreasing and $\Delta_{\X}(x_{n},\epsilon_{n})>0$ for $n\in\N$, where 
$\epsilon_{n}\stackrel{\cdot}{=}2\diam(S(\B_{\X},g_{n},\alpha_{n}+\frac{1-\alpha_{n}}{2}))\rightarrow 0$ as $n\rightarrow\infty$. 
By \eqref{eq: limxy} we get that 
$\lim_{t\rightarrow \epsilon_{n}^{-}}\Delta_{\X}(x,t)=\lim_{t\rightarrow \epsilon_{n}^{-}}\Delta_{\X}(y,t)>0$ 
for any $x,y\in\S_{\X},\ n\in\N$. Hence we conclude that 
$\inf_{x\in\S_{\X}}\Delta_{\X}(x,t)>0$ for all $t>0$. It follows easily that $\X$ is uniformly convex.

Observe that if $\X$ is uniformly convex then $\X$ and $\X^{\ast}$ are reflexive. In particular, $\X^{\ast}$ has the Radon-Nikodym Property
(RNP), so that $\B_{\X^{\ast}}$ is dentable (and actually $\omega^{\ast}$-dentable, since $\S_{\X}=\S_{\X^{\ast\ast}}$).
We refer to \cite[Ch.1]{JL} for more information about the well-known geometric concepts and results used here.

Consider next the case where $\B_{\X^{\ast}}$ is $\omega^{\ast}$-dentable. We claim that $\X^{\ast}$ is asymptotically transitive.
Let $(z_{n})\subset \S_{\X}\subset\S_{\X^{\ast\ast}}$ be a sequence such that\\ 
$\diam\left(S\left(\B_{\X^{\ast}},z_{n},\frac{1}{1+2^{-n}}\right)\right)\rightarrow 0$ as $n\rightarrow\infty$. 
Suppose that $f,g\in \S_{\X^{\ast}}$ are norm-attaining functionals and $x,y\in\S_{\X}$ are such that $f(x)=g(y)=1$. 
Since $\X$ is asymptotically transitive, there exist sequences $(T_{n}),(S_{n})\subset \mathrm{Aut}(\X)$ such that 
$T^{-1}_{n}(x)=S^{-1}_{n}(y)=z_{n}$ and $\max(||T_{n}||,||T_{n}^{-1}||,||S_{n}||,||S_{n}^{-1}||)<1+2^{-n}$ for $n\in\N$.
Since 
\[\frac{T_{n}^{\ast}f}{||T_{n}^{\ast}f||},\frac{S_{n}^{\ast}f}{||S_{n}^{\ast}f||}\in S\left(\B_{\X^{\ast}},z_{n},\frac{1}{1+2^{-n}}\right),\quad n\in\N,\]
we obtain that $||T_{n}^{\ast}f-S_{n}^{\ast}g ||\rightarrow 0$ as $n\rightarrow\infty$. 
Hence $||(T_{n}^{\ast})^{-1}S_{n}^{\ast}g -f||\rightarrow 0$ and 
$\max(||(S_{n}^{\ast})^{-1}T_{n}^{\ast}||,||(T_{n}^{\ast})^{-1}S_{n}^{\ast}||)\rightarrow 1$ by assumption as $n\rightarrow\infty$.
By the Bishop-Phelps Theorem the norm-attaining functionals are dense in $\S_{\X^{\ast}}$, and since $\mathcal{O}(h)$ is norm closed for 
all $h\in\S_{\X^{\ast}}$, we get that $\X^{\ast}$ is asymptotically transitive.
Thus an application of the first part of the proof gives that $\X^{\ast}$ is uniformly convex.

We conclude that in both our cases $\X$ and $\X^{\ast}$ are reflexive. 
It follows by applying the RNP together with the previous arguments that in fact both $\X$ and $\X^{\ast}$ 
are asymptotically transitive and uniformly convex. Hence $\X$ and $\X^{\ast}$ are also uniformly smooth.
\end{proof}

\section{Projections onto subspaces of $L^{p}$}
In this section we obtain, perhaps surprisingly, some information
about the classical $L^{p}$ spaces that appears to be new. Note that the $L^{p}$ spaces are in particular asymptotically 
transitive as they are almost transitive (see e.g. \cite[p. 8]{GR}). Hence the main idea here is to apply the rotational structure of 
$L^{p}$ to study its other types of structures.

Let us fix some notations and recall some results, which are applied in this section. For $1<p<\infty$ we denote 
$L^{p}=L^{p}(0,1)$, by $\boldsymbol{1}\in L^{p}\cup L^{p^{\ast}}$ the unit function. 
Write $M^{p}=\ker(\boldsymbol{1})\subset L^{p}$, where we consider $\boldsymbol{1}\in L^{p^{\ast}}$.
Similarly we consider $2^{\frac{1}{p^{\ast}}}\chi_{[0,2^{-1}]}$ as a functional in $L^{p^{\ast}}$.
Recall Lamperti's result that in $L^{p}$ it holds that 
\begin{equation}\label{eq: Lam}
||x+y||_{p}^{p}+||x-y||_{p}^{p}=2(||x||_{p}^{p}+||y||_{p}^{p})
\end{equation}
if and only if the essential supports $\mathrm{supp}(x)$ and $\mathrm{supp}(y)$ are disjoint (see e.g. \cite[p.163]{Lac}). 
Recall also that a continuous linear projection \mbox{$P\colon \X\rightarrow \Y,\ $}$\Y\subset\X,$ which satisfies
$||P||=\inf\{||Q||:\ Q\colon \X \rightarrow\Y\ \mathrm{is\ a\ linear\ projection}\}$ is called a \emph{minimal} projection.
 
\begin{theorem}
For $1\leq p\leq \infty$ denote 
\[\alpha_{p}=\sup_{t\in (0,1)}(t^{p-1}+(1-t)^{p-1})^{\frac{1}{p}}(t^{p^{\ast}-1}+(1-t)^{p^{\ast}-1})^{\frac{1}{p^{\ast}}}-1,\quad 
\frac{1}{p}+\frac{1}{p^{\ast}}=1.\] 
Suppose that $P\colon L^{p}\rightarrow [x],\ x\in \S_{L^{p}},$ is a linear projection. Then the following conditions are equivalent:
\begin{enumerate}
\item[(1)]{$||P||=1$.}
\item[(2)]{$||\I-P||=1+\alpha_{p}$.}
\item[(3)]{$\I-P\colon L^{p}\rightarrow \ker(P)$ is a minimal projection.}
\end{enumerate}
\end{theorem}
\begin{proof}
The case $p=2$ is clear. Let us consider the case $1<p<\infty,\ p\neq 2$. 
By applying Proposition \ref{IP} we have that if $||P||=1$ then $||\I-P||=||I-\1\otimes \1||$.
On the other hand, by \cite[Thm.6]{RoPr} every linear projection $Q$ onto a $1$-codimensional subspace $Z\subset L^{p}$ satisfies
$||Q||\geq ||\I-\1\otimes \1||$. The exact value 
\[||\I-\1\otimes \1||=\max_{t\in [0,1]}(t^{p-1}+(1-t)^{p-1})^{\frac{1}{p}} (t^{p^{\ast}-1}+(1-t)^{p^{\ast}-1})^{\frac{1}{p^{\ast}}}>1\]
was calculated by Franchetti in \cite[Thm. 3]{Fr}. Observe that by the Riesz Lemma, the Hahn-Banach Theorem and the weak compactness of 
$\B_{L^{p}}$ there is for each $1$-codimensional subspace $Z\subset L^{p}$ a point $y\in S_{L^{p}}$ and a linear $1$-dimensional 
norm-one projection $P\colon L^{p}\rightarrow [y]$ such that $\ker(P)=Z$. 
Finally, recall that in a uniformly convex space each 
\textit{minimal} linear projection $Q,\ ||Q||>1,$ onto a $1$-codimensional subspace is in fact \emph{unique}, see \cite[p.28]{OL}.

The case $p\in\{1,\infty\}$ follows directly e.g. from the fact that $L^{1}$ and $L^{\infty}$ have the Daugavet property:
Whenever $T\colon L^{p}(0,1)\rightarrow L^{p}(0,1),\ p\in\{1,\infty\},$ is a compact operator, then it holds that $||\I+T||=1+||T||$
(see e.g. \cite[p. 78]{We}). On the other hand $\alpha_{p}=1$ for $p\in\{1,\infty\}$. Observe that any non-trivial projection $P$ satisfies
$||P||\geq 1$. Hence $||\I-P||\geq 2$ and equality holds if and only if $\I-P$ is minimal.
\end{proof}

The following result which is obtained by applying some classical facts about rotations of $L^{p}$, has also some nice 
consequences (see Theorem \ref{ZTh} below).
\begin{proposition}\label{isometricprop}
Let $Z_{1},Z_{2}\subset L^{p}$ be $1$-codimensional subspaces, where\\ 
$1<p<\infty$. Then $Z_{1}\oplus_{p} L^{p}=Z_{2}\oplus_{p} L^{p}$ isometrically. Moreover, \emph{exactly} one of the following holds
isometrically:
\begin{enumerate}
\item[(1)]{$Z_{1}=Z_{2}$.}
\item[(2)]{Either (2a) $Z_{1}=Z_{1}\oplus_{p} L^{p}$ or (2b) $Z_{2}=Z_{2}\oplus_{p} L^{p}$ (but not both).}
\end{enumerate}
\end{proposition}
Observe that the second condition above does occur, see Example \ref{example ni} below. 
\begin{proof}
The case $p=2$ is clear, so let us consider the case $p\neq 2$. 
Since $L^{p}$ is reflexive, there exists by a standard application of the Riesz lemma and the weak compactness of 
$\B_{L^{p}}$ points $x,y\in \S_{L^{p}}$ such that $\dist(x,Z_{1})=\dist(y,Z_{2})=1$. By the Hahn-Banach theorem there are 
$f,g\in\S_{L^{p^{\ast}}}$ such that $\ker(f)=Z_{1},\ \ker(g)=Z_{2}$ and $f(x)=g(y)=1$.

We apply the fact that there are \emph{exactly} two disjoint orbits in $\S_{L^{p}}$:
\begin{equation}\label{eq: Ba} 
\{x\in \S_{L^{p}}|m(\mathrm{supp}(x))=1\}\ \mathrm{and}\ \{x\in \S_{L^{p}}|m(\mathrm{supp}(x))<1\}, 
\end{equation}
(see \cite[p.178]{Ba}). Hence there exists a rotation $T\in\mathcal{G}_{L^{p}(0,2)}$ for which 
$T((x,0))=(y,0)$ under the identification $(x,0),(y,0)\in L^{p}\oplus_{p} L^{p}=L^{p}(0,2)$. 
The corresponding support functionals are $\tilde{f}=(f,0),\tilde{g}=(g,0)\in L^{p^{\ast}}\oplus_{p^{\ast}}L^{p^{\ast}}=L^{p^{\ast}}(0,2)$. 
Since $L^{p}(0,2)$ is smooth, it follows that $\tilde{f}=\tilde{g}\circ T$. Hence
\[Z_{1}\oplus_{p} L^{p}=\ker(\tilde{f})=\ker(\tilde{g}\circ T)=\ker(\tilde{g})=Z_{2}\oplus_{p} L^{p}\] 
isometrically, which is the first part of the claim. This also means that the conditions (2a) and (2b) are mutually exclusive 
in the case where $Z_{1}$ and $Z_{2}$ are non-isometric.

If $T\in\mathcal{G}_{L^{p}}$ is such that $T(x)=y$ then because of the smoothness of $L^{p}$ we have as above that $f=g\circ T$, so that
\begin{equation}\label{eq: Zker}
Z_{1}=\ker(f)=\ker(g\circ T)=\ker(g)=Z_{2}\ 
\end{equation}
isometrically. Suppose that $Z_{1}$ and $Z_{2}$ are \textit{non-isometric}. 
Then by relabelling we may assume without loss of generality that
\begin{equation}\label{eq: xsupp}
m(\mathrm{supp}(x))=1\ \mathrm{and}\ m(\mathrm{supp}(y))<1.
\end{equation}  
This means that $2^{\frac{1}{p}}\chi_{[0,\frac{1}{2}]}\in\mathcal{G}_{L^{p}}(y)$.
Hence we obtain that $\ker(g)$ is isometric to $M^{p}\oplus_{p} L^{p}$. 
Clearly $(M^{p}\oplus_{p} L^{p})\oplus_{p} L^{p}$ and $M^{p}\oplus_{p} L^{p}$ are isometric. Thus (2b) holds.
\end{proof}

Note the difference between the $L^{p}$-summands appearing in the proposition above and the $1$-complemented subspaces appearing 
in the next theorem.  

\begin{theorem}\label{ZTh}
Let $1<p<\infty$ and $Z\subset L^{p}$ be a finite codimensional subspace. Then there is a subspace 
$N\subset Z$ isometric to $L^{p}$ such that $N\subset L^{p}$ is $1$-complemented. 
\end{theorem}
\begin{proof}
We first verify the following claim.\\
\textit{The space $M^{p}$ contains an isometric copy of $L^{p}$ which is $1$-complemented in $L^{p}$.}\\
Indeed, 
\[N=\{f\in M^{p}|\ f(t)=-f(2^{-1}+t)\ \mathrm{for\ a.e.}\ t\in [0,2^{-1}]\}\] 
is such a subspace. Note that $T\colon N\rightarrow L^{p}(0,2^{-1})$ given by $T\colon f\mapsto 2^{\frac{1}{p}}f_{|[0,2^{-1}]}$
is an isometric isomorphism. 
The required projection $P\colon M^{p}\rightarrow N$ is given by
$P(f)(t)=2^{-1}(f(t)-f(2^{-1}+t))$ for $t\in [0,2^{-1}]$ and $P(f)(t)=-2^{-1}(f(t-2^{-1})-f(t))$ for $t\in [2^{-1},1]$.
Indeed, clearly $P$ is a linear projection. 
To verify that $||P||=1$ we use the obvious estimate $|a-b|\leq 2^{\frac{1}{p^{\ast}}}(|a|^{p}+|b|^{p})^{\frac{1}{p}}$
to obtain that 
\begin{eqnarray*}
||Pf||_{p}^{p}&=&\int_{0}^{2^{-1}}|P(f)(t)|^{p}\ dt+\int_{2^{-1}}^{1}|P(f)(t)|^{p}\ dt\\
          &=&2^{1-p}\int_{0}^{2^{-1}}|f(t)-f(2^{-1}+t)|^{p}\ dt\\
          &\leq& 2^{1-p}2^{\frac{p}{p^{\ast}}}\int_{0}^{2^{-1}} |f(t)|^{p}+|f(2^{-1}+t)|^{p}\ dt=||f||_{p}^{p}.
\end{eqnarray*} 
Hence $||P||=1$ and we have the Claim.

Suppose that $Z_{1}\subset Z_{2}\subset ...\subset Z_{n}=L^{p}$ are subspaces such that $\dim(Z_{i}\quotient Z_{i-1})=1$
for all $i-1,i\in \{1,...,n\}$. We proceed inductively.\\ 
\textit{Step 1.} In what follows we apply Proposition \ref{isometricprop} together with its proof.
We know that $Z_{n-1}$ is isometric to $M^{p}$ or $M^{p}\oplus_{p} L^{p}$.
Recall that according to the observations \eqref{eq: Ba} and \eqref{eq: Zker} by applying a suitable rotation 
$T\colon L^{p}\rightarrow L^{p}$ we may assume without loss of generality
that $Z_{n-1}\subset L^{p}$ \emph{is} $M^{p}$ or $\ker(2^{\frac{1}{p^{\ast}}}\chi_{[0,2^{-1}]})$.

If $Z_{n-1}=M^{p}$ the previous claim gives a linear norm-$1$ projection\\ 
$P_{n-1}\colon L^{p}\rightarrow N_{n-1}$, where $N_{n-1}\subset Z_{n-1}$ is an isometric copy of $L^{p}$. 

If $Z_{n-1}=\ker(2^{\frac{1}{p^{\ast}}}\chi_{[0,2^{-1}]})$ the operator $P_{n-1}\colon L^{p}\rightarrow \{f\in L^{p}|f_{|[0,2^{-1}]}\equiv 0\ \mathrm{a.e.}\}$ defined by $P_{n-1}(f)=\chi_{[2^{-1},1]}f$ is also a linear norm-$1$ projection onto.

Hence, in both the cases there exists a linear norm-$1$ projection\\ $P_{n-1}\colon L^{p}\rightarrow N_{n-1}$,
where $N_{n-1}$ is an isometric copy of $L^{p}$.\\  
\textit{Step 2.} Observe that $\dim(N_{n-1}\quotient N_{n-1}\cap Z_{n-2})\leq 1$. Since $N_{n-1}$ is isometric
to $L^{p}$ we may apply Step 1 to conclude that there exists a subspace $N_{n-2}\subseteq N_{n-1}\cap Z_{n-2}$,
which is isometric to $L^{p}$ and $1$-complemented in $N_{n-1}$. 
Denote the corresponding norm-$1$ projection by $P_{n-2}\colon N_{n-1}\rightarrow N_{n-2}$.

We continue in this manner to define subspaces $N_{n-1}\supset N_{n-2}\supset...\supset N_{1}$ which are isometric
to $L^{p}$ together with the norm-$1$ projections $P_{i-1}\colon N_{i}\rightarrow N_{i-1}$ for $i-1,i\in \{1,...,n\}$.
Hence $N_{1}\subset Z_{1}$ is the required subspace isometric to $L^{p}$ and 
$P_{1}\circ P_{2}\circ ...\circ P_{n-1}\colon L^{p}\rightarrow N_{1}$ is the corresponding norm-$1$ projection.
\end{proof}

Recall that $L^{p}$ is primary, that is, if $L^{p}=M\oplus N$, then either $M=L^{p}$ or $N=L^{p}$ isomorphically
(see e.g. \cite[2.d.11]{LT}). Hence we know that in the following result $\ker(P_{1})$ is \emph{isomorphic}
to $\ker(P_{2})$. The crux of the following result is that $\ker(P_{1})$ and $\ker(P_{2})$ below are \emph{almost isometric}.
Thus we obtain examples of subspaces of $L^{p}$ (hence both uniformly convex and uniformly smooth), which are mutually non-isometric 
but still almost isometric (see Example \ref{example uusi}).

\begin{theorem}\label{almisom}
Let $n\in\N,\ 1<p<\infty,\ p\neq2$ and suppose that $Y_{1},Y_{2}\subset L^{p}$ are isometric copies of $\ell_{n}^{p}$ or of $\ell^{p}$. 
Then there exist unique linear projections $P_{1}\colon L^{p}\rightarrow Y_{1},\ P_{2}\colon L^{p}\rightarrow Y_{2}$ such that
$||P_{1}||=||P_{2}||=1$ and  $\ker(P_{1})$ and $\ker(P_{2})$ are \emph{almost isometric}.
\end{theorem}  
\begin{proof}
We consider only the case where $Y_{1}$ and $Y_{2}$ are isometric copies of $\ell^{p}$, since the argument for the other cases is similar.
Denote by $S\colon \ell^{p}\rightarrow Y_{1}$ the corresponding isometry and put $\hat{e}_{k}=S(e_{k}),\ k\in\N,$ for the unit vector 
basis $(e_{k})\subset \ell^{p}$. Since $S$ is an isometry we obtain by the characterization \eqref{eq: Lam} that the vectors 
$\hat{e}_{k}\in L^{p},\ k\in\N,$ have pairwise disjoint essential supports. It is then a well-known fact that there is a
normalized sequence $(\hat{e}_{k}^{\ast})\subset L^{p^{\ast}}$, where $\mathrm{supp}(\hat{e}_{k}^{\ast})=\mathrm{supp}(\hat{e}_{k})$
for $k\in\N$, such that 
\[P(x)=\sum_{k\in\N}\hat{e}_{k}^{\ast}(x)\hat{e}_{k}\]
defines a norm-$1$ projection $L^{p}\rightarrow \overline{\span}(\{\hat{e}_{k}|k\in\N\})$.  
By the result of Beauzamy and Maurey \cite[Pf. of Prop. 5]{BM} this projection $P$ is unique because $L^{p}$ is smooth. 

Hence we may partition $[0,1]=\bigcup_{k\in\N} A_{k}$, where 
$A_{1}=[0,1]\setminus \bigcup_{k\geq 2}\mathrm{supp}(\hat{e}_{k})$,
$A_{2}=\mathrm{supp}(\hat{e}_{2})$, $A_{3}=\mathrm{supp}(\hat{e}_{3})$,...
Note that $\mathrm{supp}(\hat{e}_{k}^{\ast})\subset A_{k}$ for $k\in\N$.  
For each $\hat{e}_{k\ |A_{k}}^{\ast}\in L^{p^{\ast}}(A_{k})$ it holds that 
$\ker\left(\hat{e}_{k\ |A_{k}}^{\ast}\right)\subset L^{p}(A_{k})$ 
is a $1$-codimensional subspace.

Under the above notations one can write 
\[\ker(P)=\bigoplus_{k\in\N}\ker\left(\hat{e}_{k\ |A_{k}}^{\ast}\right)\subset L^{p}, \]
where the direct sum is understood in the $\ell^{p}$-sense. It suffices for the claim to show that 
\[d_{\mathrm{BM}}\left(\ker\left(\hat{e}_{k\ |A_{k}}^{\ast}\right),\ker\left(\1_{|A_{k}}\right)\right)=1\] 
for any $k\in\N$, where we consider $\1_{|A_{k}}\in L^{p^{\ast}}(A_{k})$. 
Indeed, if $T_{k}\colon \ker\left(\1_{|A_{k}}\right)\rightarrow \ker\left(\hat{e}_{k\ |A_{k}}^{\ast}\right),\ k\in\N,$ 
are isomorphisms such that $||T_{k}||\leq 1+C2^{-k},\ C\geq 0,$ and $||T_{k}^{-1}||=1$, then also 
\[\boldsymbol{T}=\underset{k\in\N}{\bigoplus} T_{k}\colon \underset{k\in\N}{\bigoplus}\ \ker\left(\1_{|A_{k}}\right)\rightarrow \underset{k\in\N}{\bigoplus} \ker\left(\hat{e}_{k\ A_{k}}^{\ast}\right),\]
(where the direct sums are in $\ell^{p}$-sense) satisfies that $||\boldsymbol{T}^{-1}||=1$ and $||\boldsymbol{T}||\leq 1+C$. 
Above the domain and the range of $\boldsymbol{T}$,
respectively $\underset{k\in\N}{\bigoplus}\ \ker\left(\1_{|A_{k}}\right)$ and $\ker(P)$ lie in $L^{p}$. 
Even though the estimate above is by no means sharp, it is suitable for the argument at hand.

Fix $k_{0}$. We proceed to show that $d_{\mathrm{BM}}\left(\ker\left(\hat{e}_{k_{0}\ |A_{k_{0}}}^{\ast}\right),\ker\left(\1_{|A_{k_{0}}}\right)\right)=1$.
As $L^{p}(A_{k_{0}})$ and $L^{p}(0,1)$ are isometric, we may assume without loss of generality that $A_{k_{0}}=[0,1]$. 
Since $L^{p}$ is asymptotically transitive, there is by definition a sequence $(T_{k})\subset \mathrm{Aut}(L^{p})$ of
automorphisms such that $T_{k}(\hat{e}_{k_{0}})=\1$ for all $k\in\N$ and such that $\max(||T_{k}||,||T_{k}^{-1}||)\rightarrow 1$
as $k\rightarrow\infty$. Observe that $(\I-\1\otimes \1)_{|\ker(\1)}=\I_{|\ker(\1)}$.
Define $S_{k}\in L(L^{p},\ker(\1))$ for $k\in\N$ by $S_{k}=(\I-\1\otimes \1)\circ T_{k}$. 
Since $\hat{e}_{k_{0}}^{\ast}\circ T_{k}^{-1}(\1)=1$ for $k\in\N$ and 
$||\hat{e}_{k_{0}}^{\ast}\circ T_{k}^{-1}||\rightarrow 1$ as $k\rightarrow\infty$ we obtain by the Frechet-smoothness of $L^{p}$ 
and the Smulyan lemma (see Lemma \ref{lm FS}) that 
$\1-\hat{e}_{k_{0}}^{\ast}\circ T_{k}^{-1}\stackrel{||\cdot||}{\longrightarrow}0$ as $k\rightarrow\infty$. Hence we get that
\begin{eqnarray*}
  &&T_{k}\circ \left(\I-\hat{e}_{k_{0}}^{\ast}\otimes \hat{e}_{k_{0}}\right)-S_{k}=\left(\I-\left(\hat{e}_{k_{0}}^{\ast}\circ T_{k}^{-1}\right)\otimes T_{k}(\hat{e}_{k_{0}})\right)\circ T_{k} - S_{k}\\
  &&=\left(\I-\left(\hat{e}_{k_{0}}^{\ast}\circ T_{k}^{-1}\right)\otimes \1\right)\circ T_{k} - (\I-\1\otimes \1)\circ T_{k} \stackrel{||\cdot||}{\longrightarrow}0\ \mathrm{as}\ n\rightarrow\infty. 
\end{eqnarray*}
This observation justifies the fact that for large enough $k$ it holds that
\[\codim\left(S_{k}\left(\ker\left(\hat{e}_{k_{0}}^{\ast}\right)\right)\right)=\codim\left(\left(\I-\hat{e}_{k_{0}}^{\ast}\otimes \hat{e}_{k_{0}}\right)\left(\ker\left(\hat{e}_{k_{0}}^{\ast}\right)\right)\right)=1.\]
Since $\codim(S_{k}(L^{p}))=1$ for $k\in\N$, we obtain that for sufficiently large $k$ it holds that 
$S_{k}\left(\ker\left(\hat{e}_{k_{0}}^{\ast}\right)\right)=S_{k}(L^{p})$, so that the map 
$S_{k\ |\ker\left(\hat{e}_{k_{0}}^{\ast}\right)}\colon \ker\left(\hat{e}_{k_{0}}^{\ast}\right)\rightarrow \ker(\1)$ is onto.
Observe that 
\[\max\big(||T_{k}\circ \left(\I-\hat{e}_{k_{0}}^{\ast}\otimes \hat{e}_{k_{0}}\right)_{|\ker\left(\hat{e}_{k_{0}}^{\ast}\right)}||,||(T_{k}\circ \left(\I-\hat{e}_{k_{0}}^{\ast}\otimes \hat{e}_{k_{0}}\right)_{|\ker\left(\hat{e}_{k_{0}}^{\ast}\right)})^{-1}||\big)\rightarrow 1\]
as $k\rightarrow \infty$, and on the other hand that
\[(T_{k}\circ \left(\I-\hat{e}_{k_{0}}^{\ast}\otimes \hat{e}_{k_{0}}\right))_{|\ker\left(\hat{e}_{k_{0}}^{\ast}\right)}-S_{k\ |\ker\left(\hat{e}_{k_{0}}^{\ast}\right)}\stackrel{||\cdot||}{\longrightarrow} 0\ \mathrm{as}\ k\rightarrow\infty.\]
Thus the restriction $S_{k\ |\ker\left(\hat{e}_{k_{0}}^{\ast}\right)}\colon \ker\left(\hat{e}_{k_{0}}^{\ast}\right)\rightarrow \ker(\1)$ 
is an isomorphism for sufficiently large $k$, and moreover
\[\max(||S_{k\ |\ker\left(\hat{e}_{k_{0}}^{\ast}\right)}||,||(S_{k\ |\ker\left(\hat{e}_{k_{0}}^{\ast}\right)})^{-1}||)\rightarrow 1\ \mathrm{as}\ n\rightarrow\infty.\] 
Hence $d_{\mathrm{BM}}\left(\ker\left(\hat{e}_{k_{0}}^{\ast}\right),\ker(\1)\right)=1$.
\end{proof}
\begin{example}\label{example ni}
The subspace $M^{p}\subset L^{p}$ for $1\leq p <\infty,\ p\neq 2,$ is not isometrically of the form $N\oplus_{p} L^{p}$ for any 
closed subspace $N\neq \{0\}$.
\end{example}
Assume to the contrary that $M^{p}=N\oplus_{p} K$ isometrically for some non-trivial subspaces $N,K\subset L^{p}$,
where $K$ is isometric to $L^{p}$. By using the disjointness condition \eqref{eq: Lam} we obtain that there exists
a measurable decomposition $[0,1]=A\cup B$ such that for any functions $f\in N$ and $g\in K$ it holds that
$\mathrm{supp}(f)\subset A$ and $\mathrm{supp}(g)\subset B$. Indeed, assume to the contrary that there are $f\in N$
and $g\in K$ such that $m(\mathrm{supp}(f)\cap \mathrm{supp}(g))>0$. Then \eqref{eq: Lam} gives that
\[||f+g||_{p}^{p}+||f-g||_{p}^{p}\neq 2(||f||_{p}^{p}+||g||_{p}^{p}).\] 
But this contradicts the fact that $||f\pm g||_{p}^{p}=||f||_{p}^{p}+||g||_{p}^{p}$. 

Hence $\int_{A}f\ dt=\int_{B}g\ dt=0$. This states that 
\begin{equation*}
M^{p}=\left\{ f\in L^{p}|\int_{A}f\ dt=\int_{B}f\ dt=0 \right\}.
\end{equation*}
But this is a contradiction since there exist $a,b\in (0,\infty)$ such that\\ 
$\int_{[0,1]}a\chi_{A}-b\chi_{B}\ dt=0$ and hence $a\chi_{A}-b\chi_{B}\in M^{p}$.\ \qed
\begin{example}\label{example uusi}\ 
The subspaces $\ker\left(\boldsymbol{1}\right)$ and $\ker\left(2^{\frac{1}{p^{\ast}}}\chi_{[0,2^{-1}]}\right)$ of $L^{p}$ 
are almost isometric but not isometric.
\end{example}
We may apply Theorem \ref{almisom} for $1$-dimensional subspaces to obtain that the subspaces above are almost isometric.
Observe that $\ker\left(2^{\frac{1}{p^{\ast}}}\chi_{[0,2^{-1}]}\right)=M^{p}\oplus_{p} L^{p}$ isometrically.
On the other hand $\ker\left(\boldsymbol{1}\right)$ does not have such an isometric decomposition according to
Example \ref{example ni}. Hence $\ker\left(\boldsymbol{1}\right)$ and $\ker\left(2^{\frac{1}{p^{\ast}}}\chi_{[0,2^{-1}]}\right)$ 
are not isometric.\ \qed

\section{Concluding remarks}
Next we will give an asymptotic analogue of the following known comparison principle (see \cite[p.16]{GR}):\\
\textit{For a convex-transitive space $(\X,||\cdot||)$ the condition 
$\mathcal{G}_{||\cdot||}\subset \mathcal{G}_{|||\cdot|||}$ for some equivalent norm $|||\cdot|||\sim ||\cdot||$ implies that
$|||\cdot|||=c||\cdot||$ for some constant $c>0$.}

First we introduce asymptotic anologues for the expressions '$\mathcal{G}_{||\cdot||}$' 
and '$\mathcal{G}_{||\cdot||}\subset\mathcal{G}_{|||\cdot|||}$'. For a Banach space $(\X,||\cdot||)$ we denote 
\[\mathcal{F}_{||\cdot||}(\delta)=\{T\in \aut(\X,||\cdot||):\max(||T||,||T^{-1}||)\leq 1+\delta\}\]
for $\delta\geq 0$. We denote the increasing family $\{\mathcal{F}_{||\cdot||}(\delta)\}_{\delta\geq 0}$ by $\mathcal{F}_{||\cdot||}$.
If $||\cdot||$ and $|||\cdot|||$ are norms on $\X$, then the condition that for each $\epsilon>0$ there exists
$\delta>0$ such that $\mathcal{F}_{||\cdot||}(\delta)\subset \mathcal{F}_{|||\cdot|||}(\epsilon)$ 
will be denoted by $\mathcal{F}_{||\cdot||}<<\mathcal{F}_{|||\cdot|||}$. 
\begin{proposition}
Let $(\X,||\cdot||)$ be an asymptotically transitive normed space and let $|||\cdot|||$ be \emph{any} norm on $\X$. If 
$\mathcal{F}_{||\cdot||}<<\mathcal{F}_{|||\cdot|||}$, then $|||\cdot|||=c||\cdot||$ for some $c>0$.
\end{proposition}
\begin{proof}
Let $x,y\in\S_{||\cdot||}$ be such that $|||x|||\leq |||y|||$.
The assumptions yield that for each $\epsilon>0$ there exists $\delta\in (0,\epsilon)$ such that 
$\mathcal{F}_{||\cdot||}(\delta)\subset \mathcal{F}_{|||\cdot|||}(\epsilon)$. 
The asymptotic transitivity of $(\X,||\cdot||)$ implies that there is for each $\delta>0$ an automorphism
$T\in \mathcal{F}_{||\cdot||}(\delta)$ such that $T(x)=y$. This means that 
\[|||y|||-|||x|||\leq (|||T|||-1)|||x|||\leq\epsilon|||x|||.\] 
Since $\epsilon$ was arbitrary, we obtain that
$|||x|||=|||y|||$. Since $x,y\in\S_{||\cdot||}$ were arbitrary, we conclude that $|||\cdot|||=c||\cdot||$ for some $c>0$.
\end{proof}

We would like to stress the significance of the following problem.
\begin{problem}
Does there exist an asymptotically transitive Banach space $\X$, which is not almost transitive?
\end{problem}
The following example shows that for a given $x\in\S_{\X}$ the generalized orbit $\mathcal{O}(x)$ does not necessarily coincide with 
the closure of the regular orbit $\overline{\mathcal{G}(x)}$.
\begin{example}\label{example XYZ}\ 
Let $p\in (1,\infty),\ p\neq 2,$ and $M,N\subset L^{p}$ be the $1$-codimensional spaces appearing in 
Example \ref{example uusi}. Recall that $M$ and $N$ are almost isometric and non-isometric. Put 
\[\X=(M\oplus_{p}\R)\oplus_{1}(N\oplus_{p} \R).\]
We apply notation $(f,a,g,b)\in \X$, where $f\in M$, $g\in N$ and $a,b\in\R$.
Since $M$ and $N$ are almost isometric there is for each $\epsilon>0$ an automorphism $S\colon \X\rightarrow \X$ such that 
$\max(||S||,||S^{-1}||)\leq 1+\epsilon$ and $S((0,1,0,0))=(0,0,0,1)$.
Fix $T\in \mathcal{G}_{\X}$. We claim that 
\[T((M\oplus_{p}\R)\oplus_{1} (\{0\}\oplus_{p} \{0\}))=(M\oplus_{p}\R)\oplus_{1}(\{0\}\oplus_{p} \{0\}).\]
Observe that $(f,a,0,0),(0,0,g,b)\in \S_{\X}$ are exactly all the extreme points of $\B_{\X}$.
Clearly $T$ preserves extreme points. Hence $T((f,a,0,0))$ has either the form $(f_{1},a_{1},0,0)$
or $(0,0,g_{1},b_{1})$. For any $(f_{1},a_{1}),(f_{2},a_{2})$$\subset M\oplus_{p}\R\setminus \{(0,0)\}$
such that $(f_{1},a_{1})\notin [(f_{2},a_{2})]$ it holds that 
$||(f_{1},a_{1})+(f_{2},a_{2})||_{M\oplus_{p}\R}<||(f_{1},a_{1})||_{M\oplus_{p}\R}+||(f_{2},a_{2})||_{M\oplus_{p}\R}$.
Hence $T((M\oplus_{p}\R)\oplus_{1}(\{0\}\oplus_{p} \{0\}))$ is either $(M\oplus_{p}\R)\oplus_{1}(\{0\}\oplus_{p} \{0\})$ or
$(\{0\}\oplus_{p}\{0\})\oplus_{1}(N\oplus_{p}\R)$. Observe that the set $\{\pm (0,1,0,0),\pm (0,0,0,1)\}$ can be defined in purely
metric terms. Indeed, for $z\in \{\pm (0,1,0,0),\pm (0,0,0,1)\}$ there does not exist $x,y\in \X\setminus \{0\},\ x\notin [y],$
such that $x+y=z$ and $||z||_{\X}^{r}=||x||_{\X}^{r}+||Y||_{\X}^{r}$ for any $r\in\{1,p\}$. It is easy to see that there do not exist such 
atoms in $\S_{\X}$ apart from $\{\pm (0,1,0,0),\pm (0,0,0,1)\}$. Thus, if     
$T((M\oplus_{p}\R)\oplus_{1}(\{0\}\oplus_{p} \{0\}))=(\{0\}\oplus_{p}\{0\})\oplus_{1}(N\oplus_{p}\R)$ then
$T((0,1,0,0))=\pm (0,0,0,1)$. But this gives that   
$T((M\oplus_{p}\{0\})\oplus_{1} (\{0\}\oplus_{p} \{0\}))=(\{0\}\oplus_{p}\{0\})\oplus_{1}(N\oplus_{p}\{0\})$, which
is impossible since $N$ and $M$ are non-isometric.\qed
\end{example}

The fact that $M^{p}$ and $M^{p}\oplus_{p} L^{p}$ are almost isometric leads to asking if 
$M^{p}$ could be somehow exhausted by successive almost isometric embeddings of $L^{p}$.
Hence the following problem rises.
\begin{problem}
What is $d_{\mathrm{BM}}(M^{p},L^{p})$?
\end{problem}  

Let us mention some further related open problems. The question has been raised in \cite{Wo}, \cite{Pa} and \cite{Ca} 
whether every Banach space admits an equivalent maximal norm. Another natural question of an opposite flavour is whether every 
\emph{maximally normed} Banach space, which is \emph{isomorphic} to a Hilbert space, is in fact \emph{isometric} to one. 
This is actually a stronger formulation of the following problem, which appears in \cite[p.100]{Ca2}.
\begin{problem}\label{qatii}
Suppose that $\X$ is an almost transitive Banach space which is isomorphic to a Hilbert space. 
Does it follow that $\X$ is in fact isometric to a Hilbert space?
\end{problem}
Note that if $T\colon \X\rightarrow \H$ is an isomorphism, then we can consider the version 
\mbox{$C=T(\B_{\X})\subset\H$} of $\B_{\X}$. With the aid of Theorem \ref{Beata3} the above problem reduces to the following one, 
which concerns merely the geometry of convex bodies situated in Hilbert spaces: 
\begin{problem}
Let $\H$ be a Hilbert space and $C\subset \H$ a closed convex bounded subset such that $C=-C$ and $0\in\mathrm{int}(C)$. 
Does there exist for each $\epsilon>0$ a hyperplane $A_{\epsilon}\subset\H$ and $y\in \H$ such that 
\[\dist_{H}(C,C\cap([y]+A_{\epsilon}\cap C))<\epsilon ?\]   
\end{problem}

\subsection*{Aknowledgements}
This article is part of the writer's ongoing Ph.D. work, which is supervised by H.-O. Tylli to whom I am
grateful for careful suggestions about the presentation. I am indebted to B. Randrianantoanina for suggesting the present formulation
of Theorem \ref{Beata3} in response to the author's Licentiate Thesis \cite{Lic}. I am grateful
to G. Lewicki for helpful comments regarding minimal projections. The work has been supported financially 
by the Academy of Finland projects \# 53968 and \# 12070 during the years 2003-2005 and by the Finnish Cultural Foundation in 2006.

\end{document}